\documentclass[12pt, reqno]{amsart}
\usepackage{amssymb,latexsym,amsmath,amsfonts, amsthm}
\usepackage{latexsym}
\usepackage[numbers,sort&compress]{natbib}
\usepackage{color}
\usepackage[mathscr]{eucal}
\usepackage{comment}

\usepackage[linkcolor=blue, citecolor=black]{hyperref}
\hypersetup{colorlinks=true, urlcolor=blue}
\usepackage[nameinlink]{cleveref}
\usepackage{enumerate}
\usepackage[dvipsnames]{xcolor}

\numberwithin{equation}{section}

\theoremstyle{definition}
\newtheorem{definition}{Definition}[section]

\newtheorem{remark}[definition]{Remark}
\theoremstyle{plain}
\newtheorem{theorem}[definition]{Theorem}
\newtheorem{lemma}[definition]{Lemma}
\newtheorem{proposition}[definition]{Proposition}

\newtheorem{result}[definition]{Result}

\newtheorem{corollary}[definition]{Corollary}

\newcommand{\eps}{\varepsilon}

\newcommand{\Om}{\Omega}

\newcommand{\smoo}{\mathcal{C}}
\newcommand{\hol}{\mathcal{O}}

\newcommand{\dist}{{\rm dist}}

\newcommand{\cplx}{\mathbb{C}}

\begin{document}
	\title[characterization of bounded balanced convex domain]{A characterization of bounded balanced convex domains in $\mathbb{C}^n$}
	\author{Sanjoy Chatterjee and Golam Mostafa Mondal}
	\address{Department of Mathematics and Statistics, Indian Institute of Science Education and Research Kolkata,
		Mohanpur -- 741 246}
	\email{sc16ip002@iiserkol.ac.in}
	
	\address{Department of Mathematics, Indian Institute of Science Education and Research Pune,
		Pune -- 411 008}
	
	\email{golammostafaa@gmail.com, golam.mondal@acads.iiserpune.ac.in}
	\thanks{}


	\begin{abstract}
		In this paper, we investigate the characterization of balanced bounded convex domains in $\mathbb{C}^n$ in terms of the squeezing function. As an application, we provide a characterization of the polydisc in $\mathbb{C}^n$.
	\end{abstract}

	\maketitle
	
	\section{Introduction}
	The purpose of this paper is to provide a characterization of balanced convex domains in terms of the squeezing function. The notion of the squeezing function originates in a closely related notion introduced in the works of Liu et al. \cite{LiuYauTung2004,LiuYauTung2005} and Yeung \cite{Yeung2009}. The formal definition of the squeezing function was first introduced by Deng, Guan, and Zhang \cite{DenGuZha2012}, and it is as follows: Let $\mathbb{B}^n$ denote the unit ball in $\mathbb{C}^n$, and let $\Omega\subset \mathbb{C}^n$ be a bounded domain. The squeezing function on $\Omega$  is defined as 
	\begin{align*}
		S_{\Omega}(z) :=\sup\{r:\mathbb{B}^{n}(0, r)\subset f(\Omega), f\in \hol(\Omega,\mathbb{B}^{n}), f(z)=0, f \text{ is injective}\},   
	\end{align*}
	\noindent where $\mathbb{B}^n(0, r)$ denotes the ball of radius $r$, centered at the origin and $\hol(\Omega,\mathbb{B}^{n})$ denotes the collection all holomorphic maps from $\Omega$ to $\mathbb{B}^n.$ 
	
	\smallskip
	
	In \cite{Fornaess2019}, Forn\ae ss raised the question: \textit{What is the corresponding theory of the squeezing function when the ball is replaced by the polydisc?} Addressing this question, authors in \cite{NavSan2022} have studied the squeezing function corresponding to the polydisc $\mathbb{D}^n$ in $\mathbb{C}^n.$ In \cite{RongYang2022}, the authors introduced the definition of the squeezing function corresponding to a balanced domain in $\mathbb{C}^n$. Recall that a domain $D\subset \mathbb{C}^{n}$ is a balanced domain if, for every $z \in D$ and $|\lambda| \leq 1$, we have $\lambda z \in D$. The Minkowski functional of the domain $D$ is defined as $h_{D}( 
	z):=\inf\left\{t>0:t^{-1}z\in D \right\}$. Consequently, we have $D=\{z \in \mathbb{C}^{n}: h_{D}(z)<1\}$. In the case where $D$ is a balanced convex domain, the Minkowski functional defines a norm on $\mathbb{C}^{n}$ (see  \cite[Lemma 3.2]{RongYang2022}). We denote $D(r):=\{z\in \cplx^n:h_{D}(z)<r\}$. According to \cite[Lemma 2.1]{NavSan2022}, for $0<r<1$, we have $rD=D(r)$, where $rD=\{r.z: z \in D\}$. Let $D$ be a balanced, convex domain in $\mathbb{C}^n$, and let $\Omega\subset \mathbb{C}^n$ be a bounded domain. We denote the squeezing function on $\Omega$ with respect to $D$ by $T^{D}_{\Omega},$ defined as follows: 
	\begin{align*}
		T^{D}_{\Omega}(z) :=\sup\{r:D(r)\subset f(\Omega), f\in \hol(\Omega,D), f(z)=0, f \text{ is injective}\},    
	\end{align*}
	where $\hol(\Omega,D)$ denotes the collection all holomorphic maps from $\Omega$ to $D$. When $D=\mathbb{D}^n,$ for simplicity, we will denote $T^{\mathbb{D}^n}_{\Omega}$ as $T_{\Omega}.$
	
	\medskip
	
	The squeezing function is a useful invariant for studying the complex geometry of domains; see Deng et al. \cite{DenGuZha2016} and the references therein. Among the
	works on the complex geometry of domains that utilize the squeezing function and which have appeared after \cite{DenGuZha2016}: firstly, Diederich, Fornæss, and Wold \cite{DieForWold2016}
	presented the following characterization of the unit ball in $\mathbb{C}^n$:
	\begin{result}\cite[Theorem 2.1]{DieForWold2016}\label{R:DFW}
		Let $\Omega$ be a bounded domain with $\smoo^2$-smooth boundary. Suppose there is a sequence of points $q_{j}$ approaching the boundary, so that the squeezing function $S_{\Omega}(q_j)\ge 1-\eps_{j}\dist(q_{j},\partial\Omega), \eps_{j}\to 0$ as $j\to \infty.$ Then $\Omega$ is biholomorphic to the ball. \end{result}
	
	\par Recently, Bharali, Borah, and Gorai \cite{BharaliBorahGorai2023} presented a new application of the squeezing function. They
	show how the squeezing function can be used to detect when a given bounded pseudoconvex domain in $\mathbb{C}^n,~ n \ge 2,$ is not biholomorphic to any product domain. Also, their work makes extensive use of the squeezing
	functions $T^{D}_{\Omega}$. The present paper focuses on characterizing bounded balanced convex domains in $\mathbb{C}^n$ through the utilization of the squeezing function. In particular, we investigate when a Reinhardt domain is biholomorphic to the bounded balanced convex domains in $\cplx^n$. A domain $\Omega\subset \mathbb{C}^n$ is called Reinhardt if $(z_1,\cdots,z_n)\in \Omega$ implies $(e^{i\theta_1}z_1,\cdots,e^{i\theta_n}z_n)\in \Omega$ for all $(\theta_1,\theta_2,\cdots,\theta_n)\in \mathbb{R}^n.$  A domain $D$ is called \textit{homogeneous} if the automorphism group $Aut(D)$ of $D$ acts transitively on $D$. We provide the following characterization of the bounded balanced convex domains:
	
	\begin{theorem}\label{T:Smooth_Balanced}
		Let $\Omega$ be a bounded pseudoconvex Reinhardt domain, and let $D \Subset \mathbb{C}^{n}$ be a  homogeneous balanced convex domain. Suppose there is a sequence of points $q_{j}$ approaching the boundary, so that the squeezing function $T^{D}_{\Omega}(q_j)\ge 1-\eps_{j}(\dist(q_{j},\partial\Omega))^2, \eps_{j}\to 0$ as $j\to \infty.$ Then $\Omega$ is biholomorphic to $D$.
	\end{theorem}
	
	\begin{remark}
		In \Cref{R:DFW}, it is assumed that the domain has a $\mathcal{C}^2$ smooth boundary. The smoothness of the boundary in \Cref{R:DFW} allows to apply the Kobayashi distance estimate \eqref{E:enew4} (see \cite[Lemma 2.2]{DieForWold2016}), at the boundary point, where $q_{j}$ converges . A similar estimate is available for non-smooth pseudoconvex Reinhardt domains [(2.2), \Cref{R:Kobayashi_Estimate}]. However, if the boundary point $q_0$ (say) at which the sequence $q_j$ converges is such that $q_{0}\notin (\cplx\setminus{0})^n$, the estimate [(2.2), \Cref{R:Kobayashi_Estimate}] cannot be applied directly. To address this case, we have incorporated $(\dist (q_j,b\Omega))^2$ into the hypothesis. On the other hand, if $q_{0}\in (\cplx\setminus{0})^n$, we can utilize the estimate [(2.1), \Cref{R:Kobayashi_Estimate}]. In this case, we can consider $\dist(q_{j},\partial \Omega)$ instead of $\dist(q_{j},\partial \Omega)^2.$ The proof then proceeds analogously to \Cref{T:Smooth_Balanced}, the only change needed is to use \eqref{E:enew4} of \Cref{R:Kobayashi_Estimate}, in \Cref{L:Balanced_Phi(0)_in_disc}.  Consequently, \eqref{E:Balanced_Estimate6} and \eqref{E:e4} will change similarly.
	\end{remark}

	As a corollary of \Cref{T:Smooth_Balanced}, we get the following characterization of the polydisc:
	
	\begin{corollary}\label{T:Biholo_Polydisc}
		Let $\Omega$ be a bounded pseudoconvex Reinhardt domain. Suppose there is a sequence of points $q_{j}$ approaching the boundary, so that the squeezing function $T_{\Omega}(q_j)\ge 1-\eps_{j}(\dist(q_{j},\partial\Omega))^2,~ \eps_{j}\to 0$ as $j\to \infty.$ Then $\Omega$ is biholomorphic to the polydisc.      
	\end{corollary}

	\smallskip

	\Cref{T:Smooth_Balanced} provides a class of pseudoconvex bounded Reinhardt domains that are biholomorphic to homogeneous bounded balanced convex domains. Pin\v{c}uk \cite{Pinchuk1980} showed that if $D$ and $G$ are bounded pseudoconvex domains, where $G$ has a $\mathcal{C}^{2}$-smooth boundary and $D$ has a piecewise $\mathcal{C}^{2}$-smooth boundary that is not smooth, then there exists no proper holomorphic mapping $g:G\to D$. As a consequence, $D$ and $G$ are not biholomorphically equivalent. This result of Pin\v{c}uk \cite{Pinchuk1980} is a generalization of a result by Henkin \cite{Henkin1973}, which stated that a generalized analytic polyhedron cannot be biholomorphically equivalent to a bounded domain with a $\mathcal{C}^{2}$-smooth boundary. Consequently, it is not possible to find a bounded pseudoconvex Reinhardt domain with $\mathcal{C}^2$-smooth boundary in $\mathbb{C}^n$ whose squeezing function satisfies the estimate stated in \Cref{T:Biholo_Polydisc}. This result emphasizes the importance of exploring the squeezing function with respect to polydiscs or other domains to gain deeper insights into their complex geometry and function theory. Particularly, Fornaess's question is of significant importance in this context.

	\section{Technical Results}
	In this section, we present basic notions and essential technical results that are necessary for proving the theorems in this paper. We begin with the definition of Kobayashi pseudodistance. First, we give the definition of Lempert function.
	\begin{definition}
		Let $\Omega \subset \cplx^{n}$ be a domain and $\rho: \mathbb{D} \times \mathbb{D} \to \mathbb{R}$ denotes the Poincar\'e distance on the open unit disc. The Lempert function $\delta_{\Omega} \colon \Omega \times \Omega \to [0,\infty)$ is defined by
		\begin{align*}
			\delta_{\Omega}(z,w)=\inf\left\{\rho(\zeta_{0}, \zeta_{1}): \exists \phi \in \text{Hol}(\mathbb{D}, \Omega) \text{ with } \phi(\zeta_{0})=z, ~ \phi(\zeta_{1})=w\right\}.
		\end{align*}
	\end{definition}
	The following is the definition of Kobayashi pseudodistance.
	\begin{definition}
		Let $\Omega \subset \cplx^{n}$ be a domain. The  {\em Kobayashi pseudodistance} $K_{\Omega} \colon \Omega \times \Omega \to [0,\infty)$ is the largest (pseudo) distance bounded above by the Lempert function, that is 
		\begin{align*}
			K_{\Omega}(z,w)=\inf\left\{\sum_{j=1}^{k}\delta_{\Omega}(z_{j-1},z_{j}): k \in \mathbb{N},~ z_{0}=z, z_{k}=w, z_{1},\cdots ,z_{k-1} \in \Omega\right\}
		\end{align*}
		for all $z, w \in \Omega$.
	\end{definition}
	A domain $\Om \subset \mathbb{C}^{n}$ is said to be Kobayashi hyperbolic if $(\Omega, K_{\Omega})$ forms a metric space.
	
	The infinitesimal Kobayashi metric is defined by as follows:
	\begin{definition}
		Let $\Om \subset \mathbb{C}^{n}$  be a domain. For $(z,v) \in T\Omega$ the infinitesimal Kobayashi metric is defined by
		\begin{align*}
			F_{\Omega}(z,v)=\inf\bigg\{\frac{1}{R}&: R>0, \exists\text{ a holomorphic map}~ \varphi\colon \mathbb{D}_{R} \to \Omega \text{ with}~ \varphi(0)=z,\\ &d\varphi(0)=v\bigg\},
		\end{align*}
		\noindent where $T\Omega$ denotes the tangent bundle of $\Omega$ and $\mathbb{D}_{R}:=\{z \in \cplx: |z|<R\}$.
	\end{definition}
	The following result from  \cite[Theorem 1]{Royden} allows us to compute Kobayashi distance by integrating the Kobayashi metric.
	\begin{result}
		Let $\Omega \subset \mathbb{C}^{n}$ be domain and $p,q \in \Omega$. Let for any piecewise $\smoo^{1}$ curve $\gamma\colon [0,1] \to \Omega$ with $\gamma(0)=p$, $\gamma(1)=q$, $l(\gamma)=\int_{0}^{1}F_{\Omega}(\gamma(t), \dot{\gamma}(t))\,dt$.  Then $K_{\Omega}(p,q)=\inf\{l(\gamma)| \gamma:[0,1] \to \Omega \text{~is a piecewise $\smoo^{1}$ curve with } \gamma(0)=p, \gamma(1)=q       \}$.
	\end{result}
	The following is a well-known fact.
	\begin{result}\label{R:Abate_Kobayashi_Estmte}
		Let $X, Y$ be two complex manifolds, and $f\in \hol(X,Y).$ Then $$K_{Y}(f(z),f(w)) \le K_X (z,w)$$
		for all $z, w \in X.$ Furthermore, if $X$ is a submanifold of $Y$ then $K_{Y}|X\times X \le K_X.$
		
	\end{result}

	\begin{result}\cite[Theorem 1]{Warsz2012}.\label{R:Kobayashi_Estimate}
		Let $D$ be pseudoconvex Reinhardt domain in $\cplx^n.$ Fix $z_{0}\in D$ and $\xi_{0}\in \partial D.$ Then, for some constant $C,$ the inequality 
		\begin{align}\label{E:enew3}
			K_{D}(z_{0},z)\le C-\log \dist(z,\partial D), 
		\end{align}
		holds if $z\in D$ tends to $\zeta_{0}.$ Additionally, for $\zeta_{0} \in \mathbb{C}_{*}^{n}:=(\mathbb{C}\setminus\{0\})^{n}$ the estimate can be improved to \begin{align}\label{E:enew4}
			K_{D}(z_{0},z)\le C'-\frac{1}{2}\log \dist(z,\partial D), 
		\end{align}
		where $C'$ is a constant.     
	\end{result}
	
	The following result is from Gupta and Pant \cite[Theorem 4.6]{GuptaPant2023}.
	\begin{result}\label{R:Extremum_function}
		Let $\Omega$ be a bounded domain $\cplx^n$ and $D \subset \mathbb{\cplx}^{n} $ be a balanced convex domain. Then for every $z\in \Omega$, there exists a holomorphic embedding $f:\Omega\to D$ with $f(z) = 0$ such that $D(T^{D}_{\Omega}(z))\subset f(\Omega)$.
	\end{result}
	
	Recall that the Minkowski functional of a balanced domain $D$ is defined as $h( 
	z):=\inf\left\{t>0:t^{-1}z\in D \right\}$. The next result is from \cite[Proposition, Page 250]{Kubota1994}.
	\begin{result}\label{R:Kubota_Kobayashi_Estimate}
		Let $X$ be a complex normed space and $D$ be a convex balanced domain in $X.$ Let $K_{D}$ denote the Kobayashi pseudodistance of $D.$ Then
		\begin{align*}
			\frac{1}{2}\log\frac{1+h(z)}{1-h(z)}=K_{D}(0,z),
		\end{align*}
		where $h$ is the Minkowski functional of $D.$
	\end{result}
	
	Next, we will establish a proposition that plays a crucial role in the proof of \Cref{T:Smooth_Balanced}.
	
	\begin{proposition}\label{P:Modification_Rudin}
		Let $D$ be a convex, balanced, and bounded domain in $\mathbb{C}^n$, and let $F:D\to D$ be an automorphism of $D$. For any $s\in (h(F(0)), 1)$, the following holds:
		\begin{align*}
			\frac{s-h(F(0))}{1-sh(F(0))}D\subset F^{-1}(sD),	
		\end{align*} 
		where $h:\mathbb{C}^n\to \mathbb{R}$ is the Minkowski functional of the domain $D$. 	
	\end{proposition}
	
	\begin{remark}
		The case $F(0)=0$ is available in Rudin \cite[Theorem 8.1.2]{Rudin2008}.
	\end{remark}

	\begin{proof}[Proof of \Cref{P:Modification_Rudin}]
		In view of \cite[Lemma 3.2]{RongYang2022}, it follows that if $D$ is a convex and balanced bounded domain, then the Minkowski functional $h:\mathbb{C}^n\to \mathbb{R}$ for the domain $D$ forms a norm. We  consider the norm linear space $(\mathbb{C}^{n},h)$. Let $F(0)=a\in D$, and let us define a function $\alpha:[0,1)\to \mathbb{R}$ by
		\begin{align*}
			\alpha(x)=\frac{x-h(a)}{1-xh(a)}.   
		\end{align*}
		\noindent Note that $\alpha$ is a increasing continuous function on $[0,1).$ To prove this proposition, let us fix $s\in (h(a),1)$ and $z\in \alpha(s)D.$ Clearly, there exists $t^{'}\in (0,\alpha(s))$ such that $z\in t^{'}\overline{D}.$ Then, by the intermediate value property, we conclude that there exists $t\in (h(a),s)$ such that $\alpha(t)=t^{'}.$

		\smallskip

		\noindent Now, let us consider a complex linear functional $L\colon\mathbb{C}^n\to \mathbb{C}$ with $\|L\|_{h}=1,$ where $\|L\|_{h}$ denotes the operator norm of the linear operator $L$ when considering $(\mathbb{C}^{n}, h)$ as a norm linear space, and it is defined as follows:
		\begin{align*}
			\|L\|_{h}:=\{|L(x)|: h(x)\le 1\}.
		\end{align*}
		\noindent It is worth noting that $|L(x)|<1$ for all $x$ in $D$. For $z$ and $t$ as mentioned above, let us define the function $g:\mathbb{D}\to \mathbb{D}$ as follows:
		\begin{align*}
			g(\lambda)=L(F(\lambda (\alpha(t)^{-1} z))).	
		\end{align*}
		
		\noindent  It is easy to see that $g$ is well-defined and holomorphic on $\mathbb{D}.$ By a straightforward computation, we find that
		\begin{align*}
			|g(\lambda)|\le \frac{|\lambda|+|L(a)|}{1+|L(a)||\lambda|} ~\forall \lambda\in \mathbb{D}, \forall L\in (\mathbb{C}^n)^{*} \text{ with } \|L\|_{h}=1.
		\end{align*}
		
		\noindent This implies,
		\begin{align}\label{E:e1}
			|L(F(\lambda (\alpha(t)^{-1} z)))|\le \frac{|\lambda|+|L(a)|}{1+|L(a)||\lambda|} ~\forall \lambda\in \mathbb{D}, \forall L\in (\mathbb{C}^n)^{*} \text{ with } \|L\|_{h}=1.
		\end{align}
		
		Now we have that $|L(a)|\le \|L\|_{h}\cdot h(a)\le h(a)$. Clearly, the function $\frac{\beta+x}{1+\beta x}$ is an increasing function in $x$ on $[0,1)$ for any $\beta \in [0,1) $. Therefore, we have that 
		\begin{align}\label{E:Key_Inequality}
			\frac{|\lambda|+|L(a)|}{1+|L(a)||\lambda|}&\le \frac{|\lambda|+h(a)}{1+h(a)|\lambda|}~~ \forall \lambda \in \mathbb{D}, \forall L\in (\mathbb{C}^n)^{*} \text{ with } \|L\|_{h}=1.
		\end{align}
		
		\noindent Hence, from \eqref{E:e1} and \eqref{E:Key_Inequality}, we deduce that
		
		\begin{align}\label{E:eNew1}
			|L(F(\lambda (\alpha(t)^{-1} z)))|\leq  \frac{|\lambda|+h(a)}{1+h(a)|\lambda|}~~\forall \lambda \in \mathbb{D}, \forall L\in (\mathbb{C}^n)^{*} \text{ with } \|L\|_{h}=1.
		\end{align}

		\noindent In view of  Hahn-Banach theorem, we conclude that for every $\lambda\in \mathbb{D}$ with $F(\lambda (\alpha(t)^{-1}z)))\ne 0,$ there exists $L^{z,t}_{\lambda}\in (\mathbb{C}^{n})^{*}$ with $\|L^{z,t}_{\lambda}\|_{h}=1$ such that 
		\begin{align}\label{E:Hahn_Banach}
			L^{z,t}_{\lambda}(F(\lambda (\alpha(t)^{-1} z))) =h(F(\lambda (\alpha(t)^{-1} z))).
		\end{align}

		\noindent 
		For $L = L_{\lambda}^{z,t}$, $\lambda=\alpha(t)$, and assuming $F(z) \neq 0$, we deduce from \eqref{E:Hahn_Banach} and \eqref{E:eNew1} that
		\begin{align}\label{E:Computation_hF}
			h(F(z))\le \frac{|\alpha(t)|+|h(a)|}{1+|h(a)||\alpha(t)|}.
		\end{align}

		\noindent By a simple computation, we obtain that
		\begin{align*}
			\frac{|\alpha(t)|+h(a)}{1+h(a)|\alpha(t)|}=t.
		\end{align*}
		
		\noindent Therefore, using \eqref{E:Computation_hF}, we can deduce from the last equality that
		\begin{align*}
			h(F(z))\le t.
		\end{align*}
		Note that this is also true when $F(z)=0$.
		\noindent Hence, $F(z)\in t\overline{D}\subset sD.$ Therefore, for all $s\in (h(F(0)),1)$, we establish the following:
		\begin{align*}
			\frac{s-h(F(0))}{1-sh(F(0))}D\subset F^{-1}(sD),
		\end{align*}
		which proves the result.
	\end{proof}
	
	\section{Proof of \Cref{T:Smooth_Balanced} and \Cref{T:Biholo_Polydisc}}
	We need to make some preparations to establish this theorem and accomplish this by proving the following lemmas. In the following lemmas, $\Om$, $D$, $\{q_{j}\}$, and $\{\eps_{j}\}$ are as stated in   \Cref{T:Smooth_Balanced}.
	\begin{lemma}\label{L:Balanced_Kobayashi_Estmt5}
		There exists a sequence of embeddings $\Phi_{j}\colon \Omega\to D$ such that $\Phi_{j}(q_{j})=0$, and for all $z \in \Omega$ and $j \in \mathbb{N}$, the following estimate holds:
		
		\begin{align*}
			\frac{1}{2}\log\frac{1+h(w)}{1-h(w)}\le K_{\Phi_{j}(\Omega)}(0,w), ~\forall w \in \Phi_{j}(\Omega).
		\end{align*}
		\noindent 
		In particular,
		\begin{align*}
			\frac{1}{2}\log\frac{1+h(w)}{1-h(w)}\le K_{\Phi_{j}(\Omega)}(0,w),~\forall w\in\overline{{D}(1-\eps_{j}(\dist(q_{j},\partial\Omega))^2)}, ~\forall j \in \mathbb{N}.
		\end{align*} 
		Here, $h:\mathbb{C}^n\to \mathbb{R}$ represents the Minkowski functional of the domain $D$. 
	\end{lemma}	
	
	\begin{proof}
		From \Cref{R:Extremum_function}, we obtain the existence of a sequence of embeddings $\Phi_{j}:\Omega\to D$ such that $\Phi_{j}(q_{j})=0$ and $D(T_{\Omega}^{D}(q_{j}))\subset \Phi_{j}(\Omega)$. Based on the assumption $T^{D}_{\Omega}(q_{j})>1-\eps_{j}(\dist(q_{j},\partial\Omega))^2$, it follows that:
		
		\begin{align}\label{E:setinclusion}
			\overline{D(1-\eps_{j}(\dist(q_{j},\partial\Omega))^2)}\subset D(T^{D}_{\Omega}(q_{j}))\subset \Phi_{j}(\Omega)\subset D.
		\end{align}
		
		\noindent Since for 
		$j \in \mathbb{N}$, $\Phi_{j}(\Om)$ is a submanifold of $D$, by \Cref{R:Abate_Kobayashi_Estmte}, we get that
		\begin{align*}
			K_{D}(0,\Phi_{j}(z))\le K_{\Phi_{j}(\Omega)}(0,\Phi_{j}(z))~\forall z\in \Omega.
		\end{align*}
		\noindent
		\noindent In view of  \Cref{R:Kubota_Kobayashi_Estimate}, we deduce:
		
		\begin{align*}
			\frac{1}{2}\log\frac{1+h(\Phi_{j}(z))}{1-h(\Phi_{j}(z))}\le K_{\Phi_{j}(\Omega)}(0,\Phi_{j}(z)),~\forall j \in \mathbb{N}, \forall z\in \Omega.
		\end{align*}
		
		\noindent Therefore, for every $w \in \Phi_{j}(\Omega),$ we get that \begin{align*}
			\frac{1}{2}\log\frac{1+h(w)}{1-h(w)}  \leq K_{\Phi_{j}(\Omega)}(0,w)~\forall w \in \Phi_{j}(\Omega).
		\end{align*}
		
		\noindent This, together with \eqref{E:setinclusion}, gives us the following :
		\begin{align*}
			\frac{1}{2}\log\frac{1+h(w)}{1-h(w)}\le K_{\Phi_{j}(\Omega)}(0,w),~\forall w\in\overline{{D}(1-\eps_{j}(\dist(q_{j},\partial\Omega))^2)}, ~\forall j \in \mathbb{N}.
		\end{align*}
	\end{proof}

	\begin{lemma}\label{L:Balanced_Phi(0)_in_disc}

		Let $\Phi_{j}$ be as in \Cref{L:Balanced_Kobayashi_Estmt5}. Fix $z_{0}\in \Omega$. Then, there exists $N \in \mathbb{N}$, and $c>0$ such that $\Phi_{j}(z_{0})\in {D}\left(1-\frac{(\dist(q_{j},\partial\Omega))^2}{e^{2c}}\right)$ for all $j>N$.
	\end{lemma}
	
	\begin{proof}
		
		Here $q_{j} \in \Omega$ converges to the point $q_{0}$. Hence, in view of \Cref{R:Kobayashi_Estimate}, we obtain a  $N_{0} \in \mathbb{N}$ and $c >0$ such that $K_{\Om}(z_{0},q_{j}) \leq c-\log{\dist(q_{j}, \partial \Om)}$ for all $j >N_{0}$. Since $\eps_{j} \to 0$ as $j \to \infty$, there exists  $N_{1} \in \mathbb{N}$ such that  for $j >N_{1},$ we have  $\eps_{j}<\frac{1}{e^{2c}}.$
		Therefore,
		\begin{align*}
			\overline{{D}\left(1-\frac{(\dist(q_{j},\partial\Omega))^{2}}{e^{2c}}\right)}\subset \overline{{D}(1-\eps_{j}(\dist(q_{j},\partial\Omega))^{2})}, ~\forall j>N_{1}.
		\end{align*}
		
		\noindent From \Cref{L:Balanced_Kobayashi_Estmt5}, it follows that for all $w\in\overline{{D}(1-\eps_{j}(\dist(q_{j},\partial\Omega))^{2})}$ and for all $j \in \mathbb{N},$ the following inequality holds:
		
		\begin{align}\label{E:newe1}
			\frac{1}{2}\log\frac{1+h(w)}{1-h(w)}\le K_{\Phi_{j}(\Omega)}(0,w).
		\end{align}
		
		\noindent Hence,  particularly, for $w\in \mathbb{C}^n$ with $h(w)=1-\frac{(\dist(q_{j},\partial\Omega)))^{2}}{e^{2c}}$ and $j\ge \max\{N_{0}, N_{1}\}$, we obtain the expression from \eqref{E:newe1}:
		\begin{align}\label{E:Balanced_Estimate6}
			K_{\Phi_{j}(\Omega)}(0,w)\ge  \frac{1}{2}\log\left(\frac{2-\frac{(\dist(q_{j},\partial\Omega))^{2}}{e^{2c}}}{\frac{(\dist(q_{j},\partial\Omega))^{2}}{e^{2c}}} \right).
		\end{align}
		
		\noindent Again, choose  $j>N_{2}$ such that $1>\frac{(\dist(q_{j},\partial\Omega))^2}{e^{2c}}$. Consequently, for all  $j>N:=\max\{N_{0},N_{1},N_{2}\}$, we have the following :
		\begin{align*}
			\frac{2-\frac{(\dist(q_{j},\partial\Omega))^2}{e^{2c}}}{\frac{(\dist(q_{j},\partial\Omega)^2)}{e^{2c}}}\times\frac{(\dist(q_{j},\partial\Omega))^{2}}{e^{2c}}&>1.
		\end{align*}
		\noindent This implies:
		\begin{align*}
			\frac{2-\frac{(\dist(q_{j},\partial\Omega))^2}{e^{2c}}}{\frac{\dist(q_{j},\partial\Omega)}{e^{2c}}}&>\frac{e^{2c}}{(\dist(q_{j},\partial\Omega))^{2}},
		\end{align*}
		
		\noindent which further implies:
		
		\begin{align}\label{E:e4}
			\frac{1}{2} \log\frac{2-\frac{(\dist(q_{j},\partial\Omega))^2}{e^{2c}}}{\frac{\dist(q_{j},\partial\Omega)}{e^{2c}}}&>\frac{1}{2}\log\left(\frac{e^{2c}}{(\dist(q_{j},\partial\Omega))^{2}}\right)\notag\\
			&=c-\log(\dist(q_{j},\partial\Omega))\notag\\
			&\ge K_{\Omega}(z_{0},q_{j})\notag\\
			&=K_{\Phi_{j}(\Omega)}(\Phi_{j}(z_{0}),\Phi_{j}(q_{j}))\notag\\
			&=K_{\Phi_{j}(\Omega)}(0,\Phi_{j}(z_{0})).
		\end{align}
		
		\noindent Therefore, using \eqref{E:Balanced_Estimate6} and \eqref{E:e4}, we conclude that for all $j>N$, the following holds:
		\begin{align}\label{E:e3}
			K_{\Phi_{j}(\Omega)}(0,w)> K_{\Phi_{j}(\Omega)}(0,\Phi_{j}(z_{0})), \quad \forall w\in\Phi_{j}(\Omega) \text{ with }h(w)=1-\frac{(\dist(q_{j},\partial \Omega))^{2}}{e^{2c}}.
		\end{align}
		
		\noindent If there exists 
		$j>N$ such that  $\Phi_{j}(z_{0})\notin {D}\left(1-\frac{(\dist(q_{j},\partial \Omega))^{2}}{e^{2c}}\right)$, then the following two cases can occur:
		
		\noindent{\bf Case I}: There exists $j>N$ such that 
		\begin{align*}
			|h(\Phi_{j}(z_{0}))|=1-\frac{(\dist(q_{j},\partial \Omega))^{2}}{e^{2c}}.
		\end{align*}
		\noindent	Then, taking $w=\Phi_{j}z_{0}),$ we get from equation (\ref{E:e3}) that $K_{\Phi_{j}(\Omega)}(0,\Phi_{j}(z_{0})) > K_{\Phi_{j}(\Omega)}(0,\Phi_{j}(z_{0}))$. That is impossible.
		
		\smallskip	
		
		\noindent {\bf Case: II}: There exists $j>N$ such that 
		\begin{align*}
			|h(\Phi_{j}(z_{0}))|>1-\frac{(\dist(q_{j},\partial \Omega))^{2}}{e^{2c}}.
		\end{align*}
		We consider a piecewise $\smoo^{1}$-smooth  curve $\gamma\colon [0,1]\to \Phi_{j}(\Omega)$ such that $\gamma(0)=0$ and $\gamma(1)=\Phi_{j}(z_{0}).$ Since $h(\gamma(0))=0$ and $h(\gamma(1))>1-\frac{(\dist(q_{j},\partial \Omega)^{2})}{e^{2c}},$ by intermediate value property of continuous function, there exists $t_{0}\in (0,1)$ such that $h(\gamma(t_0))=1-\frac{(\dist(q_{j},\partial \Omega))^{2}}{e^{2c}}.$ 
		
		\noindent  Let us perform some calculations here:
		\begin{align*}
			\int^{1}_{0} {F}_{\Phi_{j}(\Omega)}(\gamma(t),\dot{\gamma}(t))dt&=\int^{t_{0}}_{0} {F}_{\Phi_{j}(\Omega)}(\gamma(t),\dot{\gamma}(t))dt+\int^{1}_{t_0}{F}_{\Phi_{j}(\Omega)}(\gamma(t),\dot{\gamma}(t))dt\\
			&\ge   K_{\Phi_{j}(\Omega)}(0,\gamma(t_{0}))+\int^{1}_{t_0}{F}_{\Phi_{j}(\Omega)}(\gamma(t),\dot{\gamma}(t))dt.
		\end{align*}
		
		\noindent This, along with the use of (\ref{E:e3}), leads to:
		\begin{align*}
			\int^{1}_{0}{F}_{\Phi_{j}(\Omega)}(\gamma(t),\dot{\gamma}(t))	dt&> K_{\Phi_{j}(\Omega)}(0,\Phi_{j}(z_0))+\int^{1}_{t_0} {F}_{\Phi_{j}(\Omega)}(\gamma(t),\dot{\gamma}(t))dt\\
			&>  K_{\Phi_{j}(\Omega)}(0,\Phi_{j}({z_{0}}))+ K_{\Phi_{j}(\Omega)}(\gamma(t_{0}),\Phi_{j}({z_{0}})).
		\end{align*}
		
		\noindent Let $A=\partial D\left(1-\frac{(\dist(q_{j},\partial \Omega)}{e^{2c}}\right):=\left\{w \in \mathbb{C}^{n}: h(w)=1-\frac{(\dist(q_{j},\partial \Omega)}{e^{2c}}\right\}$. Here, $\Phi_{j}(\Om)$ is Kobayashi hyperbolic. Hence, from \cite[Theorem, Page 1]{B72}, it follows that the distance $K_{\Phi_{j}(\Omega)}$ induces the Euclidean topology on $\Phi_{j}(\Omega)$. Consequently, $A$ is a compact subset in the metric space $(\Phi_{j}(\Omega), K_{\Phi_{j}(\Omega)})$.   Since $\Phi_{j}({z_{0}}) \notin A$ and  $A$ is a compact subset of $\Phi_{j}(\Omega)$, we have that 
		\begin{align}\label{E:newe2}
			\alpha=\displaystyle\inf_{z \in A} K_{\Phi_{j}(\Omega)}(z,\Phi_{j}({z_{0}}))>0. 
		\end{align}
		
		\noindent From the above computations, we obtain that
		\begin{align}\label{E:kobaya1}
			\int^{1}_{0} {F}_{\Phi_{j}(\Omega)}(\gamma(t),\dot{\gamma}(t))dt
			&>  K_{\Phi_{j}(\Omega)}(0,\Phi_{j}({z_{0}}))+ \alpha.
		\end{align}
		\noindent Since, in the above equation, $\gamma$ is an arbitrary piecewise  $\smoo^{1}$ curve joining the points origin and $\Phi_{j}({z_{0}})$, hence, from the definition of the Kobayashi distance between two points and \eqref{E:kobaya1}, it follows that 
		\begin{align*}
			K_{\Phi_{j}(\Omega)}(0,\Phi_{j}({z_{0}}))
			&\geq  K_{\Phi_{j}(\Omega)}(0,\Phi_{j}({z_{0}}))+ \alpha.
		\end{align*}
		\noindent Consequently, we have that $0 \geq \alpha$. From \eqref{E:newe2}, we conclude that this is not possible.
		
		\smallskip
		
		\noindent Combining the above two cases, we obtain that $h(\Phi_{j}({z_{0}}))<1-\frac{(\dist(q_{j},\partial \Omega))^{2}}{e^{2c}}$. Consequently, we have $\Phi_{j}(z_{0})\in D\left(1-\frac{(\dist(q_{j},\partial \Omega))^{2}}{e^{2c}}\right)$ for $j>N$.
	\end{proof}
	
	\begin{lemma}\label{L:Balanced_Kobayasi_Estimate4}
		Let $\Phi_{j}$ be as in \Cref{L:Balanced_Kobayashi_Estmt5} and $z_{0}\in \Omega$ be as \Cref{L:Balanced_Phi(0)_in_disc}. Then there exists a sequence of $\Psi_{j}\in Aut(D)$ and $c'>0$ such that $D(1-c'\eps_{j})\subset \Psi_{j}\circ\Phi_{j}(\Omega),$ and $\Psi_{j}\circ\Phi_{j}({z_{0}})=0$ for sufficiently large $j.$ 
	\end{lemma}
	
	\begin{proof}
		Let $c>0$ and $N \in \mathbb{N}$ be as \Cref{L:Balanced_Phi(0)_in_disc}. Since the automorphism group of $D$ acts transitively on $D$, we can choose $\Psi_{j}\in Aut(D)$ such that $\Psi_{j}\circ \Phi_{j}({z_{0}})=0.$ For $F=\Psi^{-1}_{j},$ and for all $s_{j}\in (h(\Phi_{j}({z_{0}})),1)$ from \Cref{P:Modification_Rudin}, we get that  
		\begin{align}\label{E:Application_Rudin}
			\frac{s_j-h(\Phi_{j}({z_{0}}))}{1-s_jh(\Phi_{j}({z_{0}}))}D\subset \Psi_{j}(s_jD).
		\end{align} 
		
		\noindent From \Cref{L:Balanced_Phi(0)_in_disc}, we obtain
		
		\begin{align*}
			h(\Phi_{j}(z_0))<1-\frac{(\dist(q_{j},\partial \Omega)^{2})}{e^{2c}}.
		\end{align*}
		\noindent With this choice of $c,$ we can say
		\begin{align}\label{E:Estimte_bj}
			h(\Phi_{j}({z_{0}}))<1-\frac{(\dist(q_{j},\partial\Omega))^{2}}{e^{2c}}<1-\eps_{j}(\dist(q_{j},\partial\Omega))^{2} ~~\forall j>N.
		\end{align}
		
		\noindent From now on, we assume $j>N$.	
		Taking $b_{j}=h(\Phi_{j}({z_{0}}))$ from \eqref{E:Application_Rudin}, we obtain
		\begin{align}\label{E:Eqn5}
			\frac{s_{j}-b_{j}}{1-s_{j}b_{j}}D\subset \Psi_{j}(s_{j}D). 
		\end{align}
		\noindent  For $s_{j}=1-\eps_{j}(\dist(q_{j},\partial\Omega)^{2})$, from \eqref{E:setinclusion}, we get
		\begin{align}\label{E:Eqn6}
			\Psi_{j}(s_{j}D)\subset \Psi_{j}(\Phi_{j}(\Omega)).
		\end{align}
		\noindent In view of (\ref{E:Eqn5}) and (\ref{E:Eqn6}), we obtain:
		\begin{align}\label{E:SmallDis_Lies Image}
			\frac{s_{j}-b_{j}}{1-s_{j}b_{j}}D\subset \Psi_{j}(s_{j}D)\subset \Psi_{j}(\Phi_{j}(\Omega)). 
		\end{align}
		
		\noindent We compute:
		\begin{align*}
			\left|\frac{s_j-b_{j}}{1-{b_{j}s_{j}}}\right|^2&=\frac{(s_j-b_{j})^2}{ (1-{b_{j}s_{j}})^2}\\
			&=\frac{s_{j}^2+b_{j}^2-2s_jb_j+(1-{b_{j}s_{j}})^2-(1-{b_{j}s_{j}})^2}{  (1-{b_{j}s_{j}})^2}\\
			&=1+\frac{s_{j}^2+b_{j}^2-1-({b_{j}s_{j}})^2}{  (1-{b_{j}s_{j}})^2}\\
			&=1+\frac{(s_j^2-1)(1-b^2_{j})}{ (1-{b_{j}s_{j}})^2}.
		\end{align*}
		\noindent By putting $s_{j}=1-\eps_{j}(\dist(q_{j},\partial \Omega))^2$ in the last equality, we get that 
		\begin{align*}
			\left|\frac{s_j-b_{j}}{1-{b_{j}s_{j}}}\right|^2&=1+\frac{(1-b_{j}^2)( \eps^2_{j}(\dist(q_{j},\partial\Omega))^4-2\eps_{j} (\dist(q_{j},\partial\Omega))^2  )}{ (1-{b_{j}s_{j}})^2}\\
			&=1-\frac{(1-b_{j}^2) \eps_{j} (\dist(q_{j},\partial\Omega))^2(2-\eps_{j}\dist(q_{j},\partial\Omega))^2   )}{(1-{b_{j}s_{j}})^2}\\
			&\ge 1-\frac{(1-b_{j}^2)2\eps_{j} (\dist(q_{j},\partial\Omega))^2}{ (1-b_{j})^2}\\
			&\ge 1-\frac{(1+b_{j})2\eps_{j} (\dist(q_{j},\partial\Omega))^2}{ (1-b_{j})}\\
			&\ge 1-\frac{4\eps_{j} (\dist(q_{j},\partial\Omega))^2}{ (1-b_{j})}.
		\end{align*}
		\noindent Using (\ref{E:Estimte_bj}), we get
		\begin{align*}
			\left|\frac{s_j-b_{j}}{1-{b_{j}s_{j}}}\right|^2 &\ge 1-4e^{2c}\eps_{j}. 
		\end{align*}
		\noindent This implies $ \left|\frac{s_j-b_{j}}{1-{b_{j}s_{j}}}\right|>1-3e^{2c}\eps_{j}.$
		From last inequality and by using (\ref{E:SmallDis_Lies Image}), we get that $D(1-c'\eps_{j})\subset \Psi_{j}\circ\Phi_{j}(\Omega)$ for $c'=3e^{2c},~~ \forall j>N.$
	\end{proof}

	\begin{proof}[Proof of \Cref{T:Smooth_Balanced}] 
		In view of \cite[Corollary 3.1]{RongYang2022}, it is enough to show  $T_{\Omega}^{D}(z_0)=1$. Note that $\Psi_{j}\circ \Phi_{j}$ is a holomorphic embedding and there exists $z_{0} \in \Omega$ with $(\Psi_{j}\circ \Phi_{j})({z_{0}})=0$ for all $j>N$, where $N\in\mathbb{N}$ as in \Cref{L:Balanced_Phi(0)_in_disc}. Applying \Cref{L:Balanced_Kobayasi_Estimate4}, we deduce that ${D}(1-3e^{2c}\eps_{j})\subset \Psi_{j}\circ \Phi_{j}(\Omega)$. As per the definition, this implies $T^{D}_{\Omega}({z_{0}})=1$, and consequently, $\Omega$ is biholomorphic to $D$.
	\end{proof}
	
	\begin{proof}[Proof of \Cref{T:Biholo_Polydisc}] 
		Note that the automorphism group $Aut(\mathbb{D}^n)$ of the polydisc $\mathbb{D}^n$ acts transitively on $\mathbb{D}^n.$ Thus, \Cref{T:Biholo_Polydisc} fulfills all the assumptions of \Cref{T:Smooth_Balanced}. Consequently, \Cref{T:Biholo_Polydisc} can be deduced from \Cref{T:Smooth_Balanced}.   
	\end{proof}

	\noindent {\bf Acknowledgements.} We sincerely thank the anonymous reviewer for the valuable comments and insightful feedback, which significantly contributed to improving the overall quality of the manuscript. We would like to express our gratitude to Sushil Gorai for providing valuable comments on our paper. The work of the first-named author is supported by a CSIR fellowship (File No-09/921(0283)/2019-EMR-I). The work of the second-named author is supported by a research grant of SERB (Grant No. CRG/2021/005884), Dept. of Science and Technology, Govt. of India.


\end{document}